\documentclass[reqno,12pt,a4paper]{amsart}
\usepackage{amssymb,amsfonts,eucal}
\usepackage{amsmath}
\usepackage{amsxtra}
\usepackage[dvips]{graphicx}

\usepackage{enumitem}
\usepackage{epic}
\usepackage{eepicemu}

\newcommand{\kre}{\kern 0sp\discretionary{--}{--}{--}\penalty 10000\hskip 0sp\relax}

\begin{document}

\title[Minimal pairs of convex sets which share a recession cone]
{Minimal pairs of convex sets which share a recession cone}

\address{Faculty of Mathematics and Computer Science\\Adam Mickiewicz University\\
Uniwersytetu Pozna\'nskiego 4, 87\\61-614 Pozna\'n, Poland}
\author{Jerzy Grzybowski}
\email{jgrz@amu.edu.pl}

\address{Faculty of Mathematics and Computer Science\\Adam Mickiewicz University\\
Uniwersytetu Pozna\'nskiego 4, 87\\61-614 Pozna\'n, Poland}
\author{Ryszard Urba\'nski}
\email{rich@amu.edu.pl}


\begin{abstract}
Robinson introduced a quotient space of pairs 
of convex sets which share their recession cone.
In this paper minimal pairs of unbounded convex sets, i.e. minimal 
representations of elements of Robinson's spaces are investigated. 
The fact that a minimal pair having property of translation 
is reduced is proved.
In the case of pairs of two-dimensional sets
a formula for an equivalent minimal pair is given, 
a criterion of minimality of a pair of sets is presented 
and reducibility of all minimal pairs is proved. 
Shephard--Weil--Schneider's criterion for polytopal summand of a compact convex set 
is generalized to unbounded convex sets. 
An application of minimal pairs of unbounded convex sets 
to Hartman's minimal representation of dc-functions is shown. 
Examples of minimal pairs of three-dimensional sets are given.
\end{abstract}

\maketitle

\centerline{
{\sc 
1. Introduction }
}

\vspace{3mm}
\noindent
For a family $\mathcal{C}(\mathbb{R}^n)$ of all nonempty closed convex 
subsets of $\mathbb{R}^n$ the addition $A+B:=\{a+b\,|\,a\in A, b\in B\}$ 
is called a \textit{Minkowski} or \textit{vector} or \textit{algebraic sum} of these sets.
For $A,B\in \mathcal{C}(\mathbb{R}^n)$ the modified addition $A\dot{+}B:=$\,cl\,$(A+B)$ 
turns the family $\mathcal{C}(\mathbb{R}^n)$ into a commutative semigroup 
with a neutral element $\{0\}$.
Moreover, for all $A,B\in\mathcal{C}(\mathbb{R}^n)$ and all $s,t\geqslant 0$ 
we have $s(tA)=s(tA)$, $t(A\dot{+}B)=tA\dot{+}tB$, 
$(s+t)A=sA\dot{+}tA$, $1A=A$, and $0A=\{0\}$. 
A relation $(A,B)\sim(C,D):\Longleftrightarrow A\dot{+}D=B\dot{+}C$ 
is not transitive because 
in $\mathcal{C}(\mathbb{R}^n)$ 
a cancellation law $A\dot{+}B=B\dot{+}C\Longrightarrow A=C$
does not hold true. Therefore, the family 
$\mathcal{C}(\mathbb{R}^n)$ cannot be embedded into a vector space.

\vspace{3mm}
\noindent
However, the family 
$\mathcal{B}(\mathbb{R}^n)$ of all nonempty closed bounded convex 
subsets of $\mathbb{R}^n$ can be embedded into a vector space, see Minkowski \cite{hM}. 
In a case of infinitely dimensional topological vector spaces 
a semigroup of nonempty closed bounded convex sets can be embedded into  
Minkowski--R\aa dstr\"om--H\"ormander space, see R\aa dstr\"om \cite{hR}, 
H\"ormander \cite{lH}, Drewnowski \cite{lD} and Urba\'nski \cite{rU}. 

\vspace{3mm}
\noindent
Quotient classes of pairs of convex sets are elements of Minkowski--R\aa dstr\"om--H\"ormander spaces. 
Sets in a given class can be arbitrarily large. 
The best representation of such a class would be inclusion-minimal pair.  
Inclusion-minimal pairs were studied by Bauer \cite{cB}, Scholtes \cite{PSU,sS}, Pallaschke \cite{GPU,GPU2009,PU,PUB} and 
by the authors \cite{jG,GKU,GUa,GUb} in connection 
with quasidifferential calculus. Quasidifferential calculus was developed by Demyanov and Rubinov \cite{DR1} 
and studied by many authors including Zhang, Xia, Gao and Wang \cite{ZXGW}
Basaeva, Kusraev and Kutateladze \cite{BKK}, Antczak 
\cite{tA}, Abbasov \cite{mA}, Dolgopolik \cite{mD} and others. 

\vspace{3mm}
\noindent
MRH spaces and basic facts about minimal pairs of convex sets 
are presented in Section~7. An embedding of a semigroup of convex sets 
is enabled by a cancellation law which was studied for its own sake by the authors \cite{GUb} 
and recently generalized to cornets by Moln\'ar and P\'ales \cite{MP}. 

\vspace{3mm}
\noindent
Robinson \cite{sR} proved an order cancellation law 
\begin{equation*}
A+B\subset B+C\Longrightarrow A\subset C.\,\,\tag{olc}
\end{equation*}
for $A,B,C$ from a family of 
unbounded closed convex sets $\mathcal{C}_V(\mathbb{R}^n)$ sharing a common recession cone $V$. 
Here, $V$ is a closed convex cone in $\mathbb{R}^n$ and a recession cone is defined 
as recc\,$A:=\{x\in \mathbb{R}^n\,|\,x+A\subset A\}$. 
A family $\mathcal{C}_V(\mathbb{R}^n)$ with Minkowski addition is a semigroup 
by Corollary 9.1.1 in \cite{tR} and as such can be embedded into a vector space. 
In this family the closed convex cone $V$ is a neutral element, 
$A\dot{+}B=A+B$, and multiplication by 0 has to be modified by $0A:=V$ for 
$A,B\in \mathcal{C}_V(\mathbb{R}^n)$. Since a cancellation law holds true, 
the relation "$\sim$" is transitive. We put $[A,B]:=[(A,B)]_{\sim}$.

\vspace{3mm}
\noindent
{\bf Theorem 1.1.} (Robinson, \cite{sR})
{\it 
The family of quotient classes $\widetilde{\mathbb{R}^n_V}:=
\mathcal{C}_V^2(\mathbb{R}^n)/_{\sim}$ with 
with the addition \,\,$[A,B]+[C,D]:=[A+C,B+D]$ and the multiplication 
\,\,$t[A,B]:= \left\{
\begin{array}{lll}
\,[tA,tB],   & t \geqslant 0 \\ 
\,[-tB,-tA], & t < 0
\end{array}
\right.$ is a smallest vector space into which the semigroup $\mathcal{C}_V(\mathbb{R}^n)$ 
can be embedded. 
}

\vspace{3mm}
\noindent
The embedding is defined by $\mathcal{C}_V(\mathbb{R}^n)\ni A\longmapsto [A,V]\in 
\widetilde{\mathbb{R}^n_V}$. In the vector space $\widetilde{\mathbb{R}^n_V}$ the 
neutral element is $[V,V]$ and the opposite element to $[A,B]$ is $-[A,B]=[B,A]$.  

\vspace{3mm}
\noindent
If the cone $V$ is trivial, i.e. $V=\{0\}$
the family $\mathcal{C}_V(\mathbb{R}^n)$ coincides with a well studied family 
$\mathcal{B}(\mathbb{R}^n)$ of all nonempty compact convex sets, i.e. of convex bodies.

\vspace{3mm}
\noindent
Robinson's theorem was generalized for closed convex sets in a Banach space by Bielawski and Tabor \cite{BT}.

\vspace{3mm}
\noindent 
Balashov and Polovinkin in their interesting paper \cite{BP} extended to unbounded sets 
the notion of generating sets. 
In a similar manner this paper extends to unbounded sets the notion of minimal pairs of sets. 

\vspace{3mm}
\noindent
In Section 2 we present a definition and a theorem of existence of minimal pairs of sets 
from $\mathcal{C}_V(\mathbb{R}^n)$ and the property of translation of minimal pairs. 
We give properties of a kernel of minimality $B_{*}$ of a pair $(A,B)$, i.e. a set of all such 
points $x$ that a pair $(A-x,B-x)$ is minimal. We also give a number of examples.

\vspace{3mm}
\noindent
In Section 3 we prove that a minimal pair of sets is reduced if and only if it has the 
property of translation.

\vspace{3mm}
\noindent
Properties of minimal pairs of two-dimensional sets are studied in Section 4. 
We give a criterion for being a summand in Proposition 4.2, a formula for an equivalent 
minimal pair in Theorem 4.4, a criterion of minimal pair in Theorem 4.5 
and prove the reducibility of all minimal pairs in Theorem 4.6. 

\vspace{3mm}
\noindent
We generalize Shephard--Weil--Schneider's criterion, i.e. Th. 3.2.11 in \cite{rS}, to polytopal 
summands of unbounded convex sets in Theorem 5.2. We also extend Bauer's criterion 
\cite{cB} of reduced pairs of polytopes to $V$-polytopes in Theorem 5.5.

\vspace{3mm}
\noindent
In Section 6 we present an application of minimal pairs of unbounded convex sets 
to a minimal, according to Hartman \cite{pH},  representation of dc-functions. 
    
\vspace{3mm}
\noindent
We complete our paper with two appendices. In Section 7 we present selected 
facts from \cite{GPU2009} on minimal pairs of bounded convex sets used in our proofs. In 
Section 8 we present Minkowski duality between convex 
sets and sublinear functions needed in the proof of Theorem 5.2.

\vspace{3mm}
\noindent 
{\bf 2. Minimal pairs of unbounded convex sets}

\vspace{3mm}
\noindent 
Let $V$ be a closed convex cone in $\mathbb{R}^n$ and $A,B\in \mathcal{C}_V(\mathbb{R}^n)$. 
A quotient class $[A,B]$ is ordered in the following way 
$$(A_1,B_1)\prec(A_2,B_2)\Longleftrightarrow A_1\subset A_2,B_1\subset B_2.$$

\vspace{3mm}
\noindent
If a recession cone $V$ is not trivial then a pair $(A+v,B+v), v\in V\setminus \{0\}$ is 
smaller than $(A,B)$ hence 
no pair $(A,B)\in \mathcal{C}^2_V(\mathbb{R}^n)$ 
is minimal. Therefore,
we say that a pair $(A,B)$ is 0-$minimal$ if $(A,B)$ is a minimal element in a subset
$\{(C,D)\in [A,B]\,|\,0\in D\}$ of a quotient class $[A,B]$. 

\vspace{3mm}
\noindent
The definition of 0-minimality seems very natural. In the case of a semigroup 
$\mathcal{B}(\mathbb{R}^n)=\mathcal{C}_{\{0\}}(\mathbb{R}^n)$ of 
bounded closed convex sets 
the existence of a minimal pair is guaranteed by the fact that a chain of compact sets 
has a nonempty intersection. In the case of $\mathcal{C}_V(\mathbb{R}^n)$ 
an intersection of a chain of sets 
containing 0 contains the cone $V$.    

\vspace{3mm}
\noindent
Every quotient class $[A,B]\in \mathcal{C}_V^2(\mathbb{R}^n)/_{\sim}$ 
contains a 0-minimal pair. The following theorem was proved by Grzybowski 
and Przybycie\'n \cite{GP} in 
much more general, possibly infinite dimensional, case. 

\vspace{3mm}
\noindent
{\bf Theorem 2.1 (existence of a $0$-minimal pair).} {\it For every pair $(A,B)\in \mathcal{C}_V^2(\mathbb{R}^n)$ with $0\in B$ 
there exists an equivalent $0$-minimal pair $(A',B')$ such that $A'\subset A,B'\subset B$.}

\vspace{3mm}
\noindent
Unlike in the case of minimal pairs of compact convex sets a pair 
$(A,B)\in \mathcal{C}_V^2(\mathbb{R}^n)$ may be $0$-minimal and a translated pair 
$(A-x,B-x)$ may not. We call a 
set $B_*:=\{x\in B\,|\,(A-x,B-x)\textup{ is 0-minimal }\}$  
a {\it kernel of minimality} of the pair $(A,B)$.
Obviously, $B_*\subset B$. 

\vspace{3mm}
\noindent
By $L_V=V\cap(-V)$ we denote the {\it subspace of lineality} of the cone $V$. 
Let us notice that for a pair $(A,B) \in \mathcal{C}_V^2(\mathbb{R}^n)$ 
we have the following equality
\begin{equation*}
\{b\in B \mid \,\, (A-b,B-b)\prec (A,B)\} = B\cap (-V).\,\,\tag{$\ast$}
\end{equation*}
The following proposition holds true.

\vspace{3mm}
\noindent
{\bf Proposition 2.2.} {\it Let $(A,B)\in\mathcal{C}_V(\mathbb{R}^n)$. If $x\in B_*$ 
then $B\cap (x-V)=x+L_V$.}

\vspace{3mm}
\noindent
{\it Proof.} Let $b\in B\cap (-V)$, then from $(\ast)$ we have $(A-b,B-b)\prec (A,B)$. 
Assume that $0\in B_*$. Then the pair $(A,B)$ is $0$-minimal and we get $B-b=B$. 
Hence $b\in L_V$ and $L_V\subset B\cap (-V)\subset L_V$. If $x\in B_*$ then 
$(A-x,B-x)$ is $0$-minimal and $(B-x)\cap(-V)=L_V$. 
\hfill{$\Box$}

\vspace{3mm}
\noindent
Proposition 2.2 says that the kernel of minimality is contained in the subset of 
minimal elements of $B$ with respect to the preorder $\leqslant_{V}$ generated by the cone $V$.
Notice also that if the cone $V$ is nontrivial then the set $B_*$ is contained in the boundary 
of $B$. 

\vspace{3mm}
\noindent
{\bf Lemma 2.3 ($B_*$ is an extreme subset of $B$).} {\it Let $(A,B)\in\mathcal{C}_V(\mathbb{R}^n)$. If $x,y\in B$ and $(x+y)/2\in B_*$ then $x,y\in B_*$.}

\vspace{3mm}
\noindent
{\it Proof.} Denote $z=(x+y)/2$.
By Theorem 2.1 there exists a $0$-minimal pair $(A'-x,B'-x)\prec (A-x,B-x)$.
Hence $z=(x+y)/2\in (B'+B)/2$ and the pair 
$$(\frac{A'}{2}+\frac{A}{2}-z,\frac{B'}{2}+\frac{B}{2}-z)\prec(A-z,B-z).$$ 
Since the pair $(A-z,B-z)$ is $0$-minimal, we obtain $B'/2+B/2=B$. By the 
cancellation law (olc) we get $B'/2=B/2$, and $B'=B$. Then the pair  
$(A-x,B-x)$ is $0$-minimal, and $x\in B_*$.
\hfill{$\Box$}

\vspace{3mm}
\noindent
{\bf Corollary 2.4.} {\it Let $(A,B)\in\mathcal{C}_V^2(\mathbb{R}^n)$. If $E\subset B$  
is a convex extreme subset of $B$ and the relative interior 
of $E$ intersects with $B_*$ then $E\subset B_*$.}

\vspace{3mm}
\noindent
Let $A\in \mathcal{C}(\mathbb{R}^n)$, $u,u_1,...,u_k\in \mathbb{R}^n$. Let 
$A(u)$ be a support set defined by $A(u)=\{a\in A\,|\,\langle a,u \rangle 
=\max\limits_{x\in A}\langle x,u \rangle\}$. Let 
$A(u_1,...,u_k)=A(u_1,...,u_{k-1})(u_k)$ be an iterated support set.
Notice that any subset of $A$ is a convex extreme subset of $A$ if and only if 
it is an iterated support set of $A$. In particular a singleton consisting of 
an extreme point of $A$ is an extreme subset of $A$.

\vspace{3mm}
\noindent
If $A,B,C,D\in\mathcal{C}_V(\mathbb{R}^n)$ and $\bar{u}=(u_1,...,u_k)\in (\mathbb{R}^n)^k$ 
and $V(u_1)=L_V$ then the following significant facts hold true 
$$(A+B)(\bar{u})=A(\bar{u})+B(\bar{u})$$ and 
$$(A,B)\sim (C,D)\Longrightarrow (A(\bar{u}),B(\bar{u}))\sim (C(\bar{u}),D(\bar{u})).$$

\vspace{3mm}
\noindent
The following proposition shows that kernels of minimality of pairs $(A,B)$ and 
$(B,A)$ "lie on the same side", respectively, of sets $B$ and $A$. 

\vspace{3mm}
\noindent
{\bf Proposition 2.5.} {\it Let $(A,B)\in\mathcal{C}_V(\mathbb{R}^n)$, 
$\bar{u}=(u_1,...,u_k)\in (\mathbb{R}^n)^k$ and $V(u_1)=L_V$. 
If $B(\bar{u})\subset B_*$ then $A(\bar{u})\subset A_*$ where 
$A_*=\{x\in A\,|\,(B-x,A-x)\textup{ is 0-minimal }\}$.}

\vspace{3mm}
\noindent
{\it Proof.} Let $y\in A(\bar{u})$. Then by Theorem 2.1 there exists 
a $0$-minimal pair $(B'-y,A'-y)\prec (B-y,A-y)$. Since $y\in A'\subset A$ and $y\in A(\bar{u})$, 
we obtain $A'(\bar{u})\subset A(\bar{u})$. Since $(B',A')\sim (B,A)$, we 
get $B'+A=A'+B$, and $B'(\bar{u})+A(\bar{u})=A'(\bar{u})+B(\bar{u})\subset A(\bar{u})+B(\bar{u})$. 
Hence by the order law of cancellation $B'(\bar{u})\subset B(\bar{u})$. Consider 
any $x\in B'(\bar{u})$. Since $B(\bar{u})\subset B_*$, the 
pair $(A-x,B-x)$ is $0$-minimal. Moreover, $x\in B'$ and $(A'-x,B'-x)\prec (A-x,B-x)$. 
Then $B'=B$, and we have just proved that $y\in A_*$ 
\hfill{$\Box$}

\vspace{3mm}
\noindent
A pair $(A,B)$ or a class $[A,B]$ is said to have a {\it property of translation 
of $0$-minimal pairs}  
if all equivalent 0-minimal pairs in $[A,B]$ are connected by translation. 
This property of translation is distinct from a property of translation 
of minimal pairs of bounded sets. If the cone $V$ is not trivial we write just 
'property of translation' because there is no possibility of misunderstanding. 

\vspace{3mm}
\noindent
For $(A,B)\in \mathcal{C}^2_{\{0\}}(\mathbb{R}^n)$ a property of translation of $0$-minimal 
pairs follows from a property of translation of minimal pairs but not the other way around. 
All pairs of flat compact convex sets from $\mathcal{C}^2_{\{0\}}(\mathbb{R}^2)$ satisfy 
the property of translation of minimal pairs \cite{cB,jG,sS}, but Example 2.10(i) presents a number 
of equivalent 0-minimal pairs not connected by translation.   

\vspace{3mm}
\noindent 
{\bf Proposition 2.6 (characterization of a kernel of $0$-minimal pair).}
{\it Let a $0$-minimal pair $(A,B)\in\mathcal{C}^2_V(\mathbb{R}^n)$ 
have the property of translation. 
Then the following assertions hold\textup{:}
\begin{enumerate}[leftmargin=*, label=\textup{(\alph*)}]
\item The set $\{(A-x,B-x) \,|\, x\in B_*\}$ is a set of all $0$-minimal pairs of the class $[A,B]$.
\item $x\in B_*$ if and only if $B\cap (x-V)=x+L_V$.
\item $B=B_*+V$.
\end{enumerate}
}

\noindent
{\it Proof.} 
(a) Let $(C,D)\in [A,B]$ be a $0$-minimal pair, then by a property of translation $D=B-z$  for some 
$z\in \mathbb{R}^n$. Since $0\in D$ we get $z=x$ for a some $x\in B$.

\vspace{2mm}
\noindent
(b) Let $(B-x)\cap (-V)=L_V$, we have $0\in B-x$. By (a) there exist a $0$-minimal pair 
$(A-z,B-z)$ such that $B-z\subset B-x$. Hence $(B-x)-(z-x) = B-z \subset B-x$. Now, by $(\ast)$ 
applied to $(A-x,B-x)$ we get $z-x\in (B-x)\cap (-V)=L_V.$  
Hence $B-z=B-x+V-(z-x)\supset B-x$ and we get $B-z=B-x$.

\vspace{2mm}
\noindent
(c) By (a) for any $b\in B$ 
there exists $x\in B_*$ such that $B-x\subset B-b$. 
Then $B + b-x\subset B$, and $b-x\in V$.
Therefore, $b=x +(b-x)\subset B_* +V$, 
and we get $B_* +V\subset B\subset B_* +V$.
\hfill{$\Box$}

\vspace{3mm}
\noindent
{\bf Remark 2.7.} Let us notice that in case of a $0$-minimal pair $(A,B)$ not having 
the property of translation the equality $B=B_*+V$ may hold true, see the pair $(\widehat{A}_0,\widehat{B}_0)$ in Example 2.10(ii), 
or not, see the pair $(\widehat{A}_1,\widehat{B}_1)$ in Example 2.10(ii). 

\vspace{3mm}
\noindent
Obviously, any pair $(A,V)$ has property of translation. Moreover, it is a unique 
0-minimal pair in a quotient class $[A,V]$.
The following example 
gives all 0-minimal pairs in a quotient class $[V,B]$.

\vspace{3mm}
\noindent
{\bf Example 2.8.} Let $n=2$, $A=V=\{0\}\times \mathbb{R}_{+}$ be a ray and 
$B=\{(x_1,x_2)\in\mathbb{R}^2\,|\,x_2\geqslant x_1^2\}$ be an epigraph of a quadratic function. 
A pair $(A,B)$ is obviously $0$-minimal. By Proposition 2.6 a pair $(A-x,B-x)$ is 0-minimal 
if and only if 
$B\cap ((x_1,x_2)-V)=\{(x_1,x_2)\}$, where $L_V=\{(0,0)\}$.
This equality holds true exactly when $x_2=x_1^2$, $x_1\in \mathbb{R}$. 
The set $B_*$ is equal to the boundary of the set $B$. Notice that $A_*=\{(0,0)\}$. 

\vspace{3mm}
\noindent
In $\mathbb{R}^3$ there exist equivalent minimal pairs not connected by translation. 
The following example was given as Example 4.1 in \cite{GP}. In that example a pair $(C,D)$ --  
not showed here --  was incorrectly presented as $0$-minimal.

\vspace{3mm}
\noindent
{\bf Example 2.9.} 
Let $V=\{x\in\mathbb{R}^3\,|\,x_1=x_2=0, x_3\leqslant 0\}$, \\ 
$B=\textup{ conv }\{(-1,-1,0),(-1,1,-1),(1,1,0),(1,-1,-1)\}+V$, \\ 
$A=\textup{ conv }(B\cup \{(-2,0,-1),(2,0,-1)\})+V$, \\
$F=\textup{ conv }\{(0,-1,-1),(0,0,0),(0,1,-1)\}+V$ and  \\
$E=\textup{ conv }(F\cup \{(-1,-1,-2),(-1,1,-1),(1,1,-2),(1,-1,-1)\})+V$. 
In Figure 2.1 we can see upper faces of sets $A,B,E,F\in \mathcal{C}_V(\mathbb{R}^3)$, where 
$V=\{x\in\mathbb{R}^3\,|\,x_1=x_2=0, x_3\leqslant 0\}$, large dots represent the origin, 
and numbers denote the third coordinate of vertices.
It can be checked that $A+F=B+E$ and that both pairs $(A,B)$ and $(E,F)$ are $0$-minimal.

\begin{displaymath}
\begin{array}{ccccc}
\setlength{\unitlength}{.5mm}
\begin{picture}(70,55)(0,0)

\qbezier(50,20)(50,20)(50,40)
\qbezier(50,20)(50,20)(70,40)
\qbezier(50,20)(50,20)(70,20)
\qbezier(70,20)(70,20)(70,40)
\qbezier(50,40)(50,40)(70,40)
\put(70,43){\makebox(0,0){0}}
\put(50,43){\makebox(0,0){-1}}
\put(50,16){\makebox(0,0){0}}
\put(70,17){\makebox(0,0){-1}}
\put(60,50){\makebox(0,0){$B$}}
\put(60,30){\makebox(0,0){\Large\bf .}}

\qbezier(0,30)(0,30)(10,20)
\qbezier(0,30)(0,30)(10,40)
\qbezier(10,20)(10,20)(30,40)
\qbezier(10,20)(10,20)(30,20)
\qbezier(10,40)(10,40)(30,40)
\qbezier(30,40)(30,40)(40,30)
\qbezier(40,30)(40,30)(30,20)
\put(30,44){\makebox(0,0){0}}
\put(10,44){\makebox(0,0){-1}}
\put(10,16){\makebox(0,0){0}}
\put(30,17){\makebox(0,0){-1}}
\put(42,27){\makebox(0,0){-1}}
\put(-2,27){\makebox(0,0){-1}}
\put(20,50){\makebox(0,0){$A$}}
\put(20,30){\makebox(0,0){\Large\bf .}}
\end{picture}
& 
\hspace{10mm}
&
\setlength{\unitlength}{.5mm}
\begin{picture}(30,55)(0,0)
\put(30,10){\makebox(0,0){\tiny\bf .}}
\put(30,50){\makebox(0,0){\tiny\bf .}}
\put(30,30){\makebox(0,0){\Large\bf .}}
\qbezier(30,10)(30,10)(30,50)
\put(30,6){\makebox(0,0){-1}}
\put(34,32){\makebox(0,0){0}}
\put(31,53){\makebox(0,0){-1}}
\put(40,50){\makebox(0,0){$F$}}

\qbezier(10,10)(10,10)(00,20)
\qbezier(00,20)(00,20)(00,40)
\qbezier(00,40)(00,40)(10,50)
\qbezier(10,50)(10,50)(20,40)
\qbezier(20,40)(20,40)(20,20)
\qbezier(20,20)(20,20)(10,10)
\qbezier(10,10)(10,10)(10,50)
\qbezier(00,40)(00,40)(20,20)
\put(10,6){\makebox(0,0){-1}}
\put(14,32){\makebox(0,0){0}}
\put(10,53){\makebox(0,0){-1}}
\put(-4,18){\makebox(0,0){-2}}
\put(-3,42){\makebox(0,0){-1}}
\put(23,18){\makebox(0,0){-1}}
\put(24,42){\makebox(0,0){-2}}
\put(-10,50){\makebox(0,0){$E$}}
\put(10,30){\makebox(0,0){\Large\bf .}}
\end{picture}

\end{array}
\end{displaymath}

\begin{center}
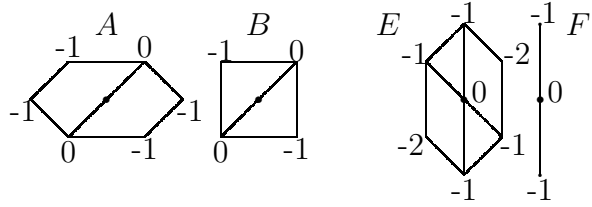
Figure 2.1. Two equivalent minimal pairs of unbounded convex sets not connected by translation from 
Example 2.9.
\end{center}

\vspace{3mm}
\noindent
Let us notice that if a given pair $(A,B)$ does not have a property of translation 
and $(B-x)\cap(-V)=L_V$ 
then a pair $(A-x,B-x)$ may or may not be 0-minimal. The following example shows 
such possibility.

\vspace{3mm}
\noindent
{\bf Example 2.10 }(i). 
Let $V=\{(0,0)\}\subset \mathbb{R}^2$, $A, B\in\mathcal{C}_V(\mathbb{R}^2)$, 
$B=$ conv $\{(0,0),(2,0)\}$ and $A=$ conv $(B\cup\{(1,1)\})$.
Let $p_0\in B$, $p_1\in \{x\in \mathbb{R}^2\,|\,x_2\leqslant \min(1-|x_1-1|,0)\}$, 
$p_2\in \{x\in \mathbb{R}^2\,|\,1-|x_1-1|< x_2< 0\}$. 
Denote $B_i=$ conv $\{\big((B-p_i)\cup \{(0,0)\}\big), i=0,1,2$, 
$A_0=A-p_0$,$A_1=$ conv $\{\big((A-p_1)\cup \{(0,0)\}\big)$, 
$A_2=$ conv $\{\big((A-p_2)\cup \{(0,0),(1,1)\}\big)$ if $(p_2)_1<0$
and $A_2=$ conv $\{\big((A-p_2)\cup \{(0,0),(-1,1)\}\big)$ if $(p_2)_1>2$. 
The sets $A_i, B_i, i=0,1,2$ are represented in Figure 2.2. 

\begin{displaymath}
\begin{array}{ccccc}
\setlength{\unitlength}{.5mm}
\begin{picture}(40,50)(0,-10)

\qbezier(00,00)(20,00)(40,00)
\qbezier(00,00)(10,10)(20,20)
\qbezier(40,00)(30,10)(20,20)
\put(20,10){\makebox(0,0){$A_0$}}
\put(00,00){\makebox(0,0){\Large\bf .}}
\put(20,-5){\makebox(0,0){$B_0$}}
\put(-2,-5){\makebox(0,0){$0$}}
\end{picture}
& 
\hspace{1mm}
&
\setlength{\unitlength}{.5mm}
\begin{picture}(40,50)(0,-10)

\qbezier(00,20)(10,30)(20,40)
\qbezier(40,20)(30,30)(20,40)
\qbezier(40,20)(20,20)(00,20)
\qbezier(30,00)(15,10)(00,20)
\qbezier(30,00)(35,10)(40,20)
\put(25,12){\makebox(0,0){$B_1$}}
\put(30,00){\makebox(0,0){\Large\bf .}}
\put(28,-5){\makebox(0,0){$0$}}
\put(45,20){\makebox(0,0){$A_1$}}
\end{picture}
& 
\hspace{1mm}
&

\setlength{\unitlength}{.5mm}
\begin{picture}(80,50)(0,-10)

\qbezier(00,00)(10,10)(20,20)
\qbezier(60,42)(40,31)(20,20)
\qbezier(60,42)(70,32)(80,22)
\qbezier(40,22)(60,22)(80,22)
\qbezier(40,22)(20,11)(00,00)
\qbezier(80,22)(40,11)(00,00)
\put(40,16){\makebox(0,0){$B_2$}}
\put(-2,-5){\makebox(0,0){$0$}}
\put(00,00){\makebox(0,0){\Large\bf .}}
\put(50,29){\makebox(0,0){$A_2$}}
\end{picture}

\end{array}
\end{displaymath}

\begin{center}
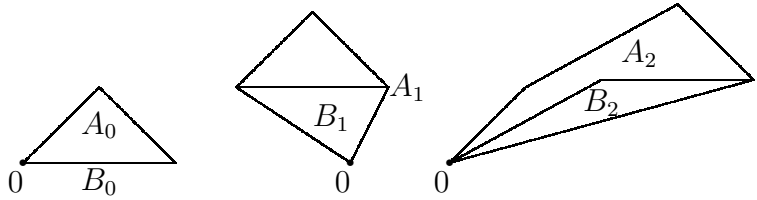
Figure 2.2. Pairs of sets described in Example 2.10(i).
\end{center}

\vspace{3mm}
\noindent
All pairs of sets $(A_i,B_i), i=0,1,2$ are equivalent to $(A,B)$. 
We are going to prove that each pair $(A_i, B_i)$ is 0-minimal. Assume that 
$(C,D)\prec (A_i, B_i)$ and $0\in D\subset B_i$.  
By Theorem 7.4 a pair of 
polygons is minimal if and only if they have at most 
one pair of parallel edges that lie on the sam side of polygons. 
Then the pair of a triangle and a segment is minimal.    
Since the segment $B$ contains $0$, the pairs $(A,B),(A_0,B_0)$ are minimal and $0$-minimal.
By Theorem 7.1, i.e. existence of equivalent minimal pair contained in a given pair,  
and by Theorem 7.2, i.e. uniqueness-up-to-translation of equivalent minimal pairs 
of flat sets,  
the set $D$ contains a translate of $B$, namely $B-p_i$. Then obviously $D=B_i$. Hence 
$A_i+B_i=A_i+D=B_i+C$,  
and by the law of cancellation $C=A_i$. Therefore, the pairs $(A_i,B_i), i=0,1,2$ are 
$0$-minimal.

\vspace{3mm}
\noindent
It can be proved that there are no other $0$-minimal pairs in the quotient class $[A,B]$. 
Notice that $(B_0)_*=B_0$, $(B_1)_*=(B_2)_*=\{(0,0)\}$. By Proposition 2.5 we obtain 
$(A_0)_*=A_0$, $(A_1)_*=\{(0,0)\}$ and $(A_2)_*=$ conv $\{(0,0),(1,1)\}$ if $(p_2)_1<0$ 
or $(A_2)_*=$ conv $\{(0,0),(-1,1)\}$ if $(p_2)_1>2$.

\vspace{3mm}
\noindent
(ii) Let $V\subset \mathbb{R}^3$ be a cone such that 
$\{x\in V\,|\,x_3\geqslant 0\}=\{(0,0,0)\}$. Denote 
$\widehat{A}=(A\times \{0\})+V$ and $\widehat{B}=(B\times \{0\})+V$. It can be proved that 
all $0$-minimal pairs in $[\widehat{A},\widehat{B}]$ are 
$(\widehat{A}_i,\widehat{B}_i), i=0,1,2$ where 
$\widehat{A}_i=(A_i\times \{0\})+V$ and $\widehat{B}_i=(B_i\times \{0\})+V$ 
and $A_i,B_i$ are sets from (i). 
Notice that $(\widehat{B}_0)_*=B_0\times \{0\}$ and $\widehat{B}_0=(\widehat{B}_0)_*+V$, 
but $(\widehat{B}_1)_*=\{(0,0,0)\}$ and $\widehat{B}_1\neq V=(\widehat{B}_0)_*+V$.

\vspace{3mm}
\noindent
{\bf 3. Reduced pairs of unbounded convex sets}

\vspace{3mm}
\noindent
Let us extend a notion of reduced pair of bounded sets from $\mathcal{B}^2(\mathbb{R}^n)$
introduced by Bauer \cite{cB}. 
A pair $(A,B)\in \mathcal{C}_V^2(\mathbb{R}^n)$ is {\it reduced} 
if $[A,B]=\{(A+M,B+M)\,|\,M\in \mathcal{C}_V(\mathbb{R}^n)\}$.

\vspace{3mm}
\noindent 
In this section we show a relationship between 
reduced pairs and the property 
of translation of 0-minimal pairs. 

\vspace{3mm}
\noindent 
{\bf Proposition 3.1.}  
{\it Let $V$ be a closed convex cone.  
If a pair $(A,B)\in \mathcal{C}_V^2(\mathbb{R}^n)$ is reduced  
then it has the property of translation. 
} 

\vspace{3mm}
\noindent 
{\it Proof.}  
Let a pair $(C,D)\in [A,B]$ be $0$-minimal. 
Then $(C,D)=(A+M,B+M)$ for some $M\in \mathcal{C}_V(\mathbb{R}^n)$. 
Since $0\in D=B+M$, there exists $b\in B$ such that $-b\in M$.
Then $A-b\subset A+M=C, B-b\subset B+M=D$. Since $(C,D)$ is $0$-minimal, 
we obtain $C=A-b, D=B-b$. 
\hfill{$\Box$}

\vspace{3mm}
\noindent 
{\bf Proposition 3.2.}   
{\it Let $V\subset \mathbb{R}^n$ be a closed convex cone.  
If a pair $(A,B)\in \mathcal{C}_V^2(\mathbb{R}^n)$ has the property of translation, 
then every  
$0$-minimal pair $(C,D)\in [A,B]$ is reduced . 
} 

\vspace{3mm}
\noindent 
{\it Proof.}  
Let a pair $(C,D)\in [A,B]$ be $0$-minimal. Let $b\in B$. 
By Proposition 2.6(a)
there exists 
$0$-minimal pair $(C-x,D-x),x\in D_*$ such that $C-x\subset A-b, D-x\subset B-b$. 
We obtain $D +b-x\subset B$, and 
$b-x\in B \dot{-}D:=\{y\,|\,D+y \subset B\}$.
Then $b=x+(b-x)\in  D+(B \dot{-}D)\subset B$, and $B\subset D+(B \dot{-}D)\subset B$. 
Hence $B=D + (B\dot{-}D)$. Then $A+D=D + (B\dot{-}D)+C$, and by the cancellation law  $A=C+(B\dot{-}D)$. 
\hfill{$\Box$}

\vspace{3mm}
\noindent 
Propositions 3.1 and 3.2 can be summed up in the following theorem.

\vspace{3mm}
\noindent 
{\bf Theorem 3.3 (equivalence of reducibility and property of translation).}  
{\it Let $V$ be a closed convex cone.  
A pair $(A,B)\in \mathcal{C}_V^2(\mathbb{R}^n)$ is reduced  
if and only if this pair has the property of translation and is a translate of some $0$-minimal pair. 
} 

\vspace{3mm}
\noindent 
{\bf Corollary 3.4.}   
{\it Let $V\subset \mathbb{R}^n$ be a closed convex cone.  
If a pair $(A,B)\in \mathcal{C}_V^2(\mathbb{R}^n)$ has the property of translation, 
then $(B,A)$ has the property of translation. 
} 

\vspace{3mm}
\noindent 
{\it Proof.}  
If $(A,B)$ has the property of translation then by Proposition 3.2 some equivalent pair 
$(C,D)$ is reduced, i.e. $[A,B]=\{(C+M,D+M)\,|\,M\in \mathcal{C}_V(\mathbb{R}^n)\}$. Hence
the pair $(D,C)$ is reduced and by Proposition 3.1 the class $[B,A]$ has the property of translation. 
\hfill{$\Box$}

\vspace{3mm}
\noindent 
Theorem 3.3 shows a difference between the property of translation of minimal pairs and 
the property of translation of $0$-minimal pairs. There exists a broad 
class of minimal pairs of compact convex sets that satisfy the property of translation not being 
reduced pairs. For example all pairs of convex polygons $(A,B)$ with exactly 
one pair $(A(u),B(u))$ of parallel edges (see Theorem 7.4). 
The authors can prove that if in a pair $(A,B)\in \mathcal{B}^2(\mathbb{R}^3)$ a set $A$ is a tetrahedron and   
for all triangular faces $A(u)$ the pairs $(A(u),B(u))$ are  
minimal then $(A,B)$ is minimal and has the property of translation.  
Such pair $(A,B)$ may be a pair of convex polyhedra and possess one or more pairs $(A(v),B(v))$ of parallel edges. 
By Theorem 5.5, i.e. Bauer's criterion for reduced polytopes the pair $(A,B)$ is not reduced.
Sufficient and necessary condition for having the property of translation in $\mathbb{R}^3$ 
is not known. On the other hand every $0$-minimal pair having the property of translation of 0-minimal pairs is 
inevitably reduced.

\vspace{3mm}
\noindent
{\bf 4. Minimal pairs of unbounded planar convex sets}

\vspace{3mm}
\noindent
In order to prove propositions and theorems of this section we need to 
present the notion and properties of an arc-length function $f_A$ corresponding to 
a planar set $A$. 
Let us consider a nonempty unbounded closed convex set $A\subset \mathbb{R}^2$. 
Let a recession cone $V$ of $A$ be pointed and unbounded . 
Obviously, $V$ is a planar convex angle of a measure $\pi-2\vartheta$ with $\vartheta \in (0,\pi/2]$. 
Assume that the negative part of the $x$-axis bisects the angle $V$.
We construct an arc-length function $f_A$ following the approach from \cite{GP}.

\vspace{3mm}
\noindent
Let $u\in \mathbb{R}^{2}$ and a support set of $A$ in the direction of $u$ be a set  
$A(u):=\{a\in A| \langle u,a \rangle =\max _{b\in A}\langle u,b \rangle\}$. 
Obviously, the support set $A(u)$ is a singleton, a segment, a ray or an empty set. 
Let $H_{A}:(-\vartheta,\vartheta)\longrightarrow \textup{bd}\,A$ be a boundary function, where 
$H_{A}(t)$ is the center of the set $A(\cos t,\sin t)$, which is either a segment or a singleton.
We also denote by 
$f_{A}:(-\vartheta,\vartheta)\longrightarrow \mathbb{R}$ 
an arc length function of $A$, with a value $f_{A}(t),t \geqslant 0$ equal to 
the length of the arc contained in the boundary $\textup{bd}\,A$ joining points
$H_{A}(0)$ and $H_{A}(t)$. If $t<0$ let a value $f_{A}(t)$
be opposite to the length of the arc joining
$H_{A}(0)$ and $H_{A}(t)$. 
The function $f_{A}$ is non-decreasing, $f_{A}(0)=0$ 
and $f_{A}(t)=\frac{1}{2}(f_{A}(t^{+})+f_{A}(t^{-}))$ 
where $f(t^{+})=\lim _{s\rightarrow t^{+}}f(s), f(t^{-})=\lim _{s\rightarrow t^{-}}f(s)$,
for $t\in (-\vartheta,\vartheta)$.

\vspace{3mm}
\noindent
On the other hand, let $f$ be any non-decreasing real function defined on an open interval 
$(-\vartheta,\vartheta)$, 
such that 
\begin{equation*}
f(0)=0 \textup{ and } f(t)=\frac{1}{2}(f(t^{+})+f(t^{-})),t\in (-\vartheta,\vartheta).\,\,\tag{$\ast\ast$}
\end{equation*}
We define the function $H_{f}:(-\vartheta,\vartheta)\longrightarrow \mathbb{R}^{2}$ 
with the help of Stieltjes integral
\begin{displaymath}
H_{f}(t):=\left\{
\begin{array}{rl}
\int\limits^{t}_{0}(-\sin s, \cos s)df(s),  & t\geqslant 0,\\
-\int\limits^{0}_{t}(-\sin s, \cos s)df(s), & t< 0.
\end{array}
\right.
\end{displaymath}
We denote $A_{f}:={\rm cl\,conv}({\rm im}H_{f})+V$.
Then we have $A_{f_{A}}=A-H_{A}(0)$ and $f_{A_{f}}=f$.
The following proposition summarizes properties 
of the correspondence between non-decreasing functions and convex sets.

\vspace{3mm}
\noindent 
{\bf Proposition 4.1.} {\it
Let $A,B\in\mathcal{C}_V(\mathbb{R}^2)$ and $f,g$ be non-decreasing functions satisfying $(\ast\ast)$.
The following formulas hold true\textup{:}\\ 
$f_{A+B}=f_A+f_B$, $f_{tA}=tf_A$ for $t\geqslant 0$,\\
$A_{f+g}=A_f+A_g$, $A_{tf}=tA_f$ for $t\geqslant 0$, \\
$f_{A_g}=g$, $A_{f_B}=B-H_B(0)=B-\textup{midpoint}B(u), u=(1,0)$,\\
$f_V\equiv 0, A_f=V$ for $f\equiv 0$.}

\vspace{3mm}
\noindent 
{\bf Proposition 4.2 (criterion of planar summands).} {\it
A set $A\in \mathcal{C}_V(\mathbb{R}^2)$ is a summand of 
$B\in \mathcal{C}_V(\mathbb{R}^2)$ if and only if a 
function $f_B-f_A$ is non-decreasing.} 

\vspace{3mm}
\noindent 
{\it Proof.} If $B=A+C$ then $f_B-f_A=f_C$.
On the other hand
if a function $g=f_B-f_A$ is non-decreasing then $B-H_B(0)=A-H_A(0)+A_g$. 
\hfill{$\Box$}

\vspace{3mm}
\noindent 
The following proposition is needed in the proof of  Theorem 5.2, i.e. a criterion for polytopal 
summands. 

\vspace{3mm}
\noindent
{\bf Proposition 4.3 (criterion of a polygonal summand).}  
{\it Let $P$ be a convex polygon, $K\in \mathcal{C}(\mathbb{R}^2)$. Assume that 
the recession cone of $K$ is not a straight line. 
Then $P$ is a summand of $K$ if and only 
if for all $u\in S^{1}$
the support set $K(u)$ is empty or contains a translate of $P(u)$.} 

\vspace{3mm}
\noindent
{\it Proof.} $\Longrightarrow )$ If $K=P+L$ for some closed convex set $L$, 
and a face $K(u)$ is nonempty 
then $K(u)=P(u)+L(u)$, and $K(u)$ contains a translate of $P(u)$.

\vspace{3mm}
\noindent
$\Longleftarrow )$ 
If the cone $V=$recc\,$K=K\dot{-}K$ is a plane or a half-plane than the theorem obviously holds true. 
Otherwise $V$ is an angle of a measure $\pi-2\vartheta$, $0<\vartheta\leqslant \pi/2$. 
We may assume that the $x$-axis is bisecting the cone $V$ and that the negative part 
of the $x$-axis is contained in $V$. 

\vspace{3mm}
\noindent
It is enough to show that $P+V$ is a summand of $K$. 
Arc-length function $f_{P+V}$ is locally constant and noncontinuous 
only at $t\in(-\vartheta,\vartheta)$ such that a support set $(P+V)(u), u=(\cos t,\sin t)$ is 
a side of $P$. Since every segment $(P+V)(u)$, having length 
equal to $f_{P+V}^{+}(t)-f_{P+V}^{-}(t)$, is 
contained in some translate of a segment $K(u)$, of the lengh 
equal to $f_K^{+}(t)-f_K^{-}(t)$, 
the difference of arc-length functions $g=f_K -f_{P+V}$ is non-decreasing.  
By Proposition 4.2 the set $P+V$ is a summand of $K$.
\hfill{$\Box$}

\vspace{3mm}
\noindent
Let us define the ordering of non-decreasing functions taking value $0$ at $0$. 
For two functions $f,g$ we say that $f$ {\it precedes} $g$ if and only if 
$g-f$ is nondecreasing. Next two theorems on $0$-minimal pairs 
in a plane correspond to Theorem 3.1 and Corollary 3.2 from \cite{GP}.

\vspace{3mm}
\noindent
{\bf Theorem 4.4 (formula for an equivalent $0$-minimal pair).}  
{\it Let $(A, B)\in \mathcal{C}_V^2(\mathbb{R}^2)$. Denote $g_A := f_A - \inf(f_A, f_B)$ and 
$g_B := f_B - \inf(f_A, f_B)$. Then the pair $(A_{g_A}+H_A(0) - H_B(0), A_{g_B})$ is 
$0$-minimal and belongs to $[A, B]$.}

\vspace{3mm}
\noindent
{\bf Theorem 4.5 (criterion of 0-minimality).}  
{\it Let $(A, B)\in \mathcal{C}_V^2(\mathbb{R}^2)$. The pair $(A,B)$ is minimal 
if and only if $\inf(f_A, f_B) \equiv 0$ and $0 \in B(cos t, sin t)$ for some 
$t \in (-\theta, \theta)$.}

\vspace{3mm}
\noindent
{\bf Theorem 4.6 ($0$-minimal pair is reduced).}  
{\it Let $V$ be a pointed unbounded convex cone in $\mathbb{R}^2$. 
Then every $0$-minimal pair $(A,B)\in \mathcal{C}_V^2(\mathbb{R}^2)$ is reduced.} 

\vspace{3mm}
\noindent 
{\it Proof.} Let $(C,D)\in [A,B]$. Then $A+D=B+C$, $f_A+f_D=f_B+f_C$ 
and $H_A(0)+H_D(0)=H_B(0)+H_C(0)$. 
We have $f_C+\inf(f_A,f_B)\prec \inf(f_C+f_A,f_C+f_B)$ and 
$\inf(f_C+f_A,f_C+f_B)-f_C \prec \inf(f_A,f_B)$.
Then $f_C+\inf(f_A,f_B)=\inf(f_C+f_A,f_C+f_B)=\inf(f_C+f_A,f_D+f_A)=f_A+\inf(f_C,f_D)$. 
Hence $g_C:=f_C-\inf(f_C,f_D)=
f_A-\inf(f_A,f_B)=f_A$. 
In a similar way $g_D:=f_D-\inf(f_C,f_D)=f_B-\inf(f_A,f_B)=f_B$. 
Thus $(C,D)=(A_{f_C}+ H_C(0), A_{f_D}+H_D(0))
=(A_{g_C}+A_{\inf(f_C,f_D)}+H_C(0), A_{g_D}+A_{\inf(f_C,f_D)}+H_D(0))=
(A_{f_A}+A_{\inf(f_C,f_D)}+H_C(0), A_{f_B}+A_{\inf(f_C,f_D)}+H_D(0))=
(A-H_A(0)+A_{\inf(f_C,f_D)}+H_C(0),B-H_B(0) +A_{\inf(f_C,f_D)}+H_D(0))=
(A+A_{\inf(f_C,f_D)}+H_D(0)-H_B(0),B+A_{\inf(f_C,f_D)}+H_D(0)-H_B(0))$. 
\hfill{$\Box$}

\vspace{3mm}
\noindent
{\bf  5. Criterion for polytopal summands}
 
\vspace{3mm}
\noindent
In this section we generalize Shephard--Weil--Schneider criterion for a polytope being 
a summand of compact convex subset of $\mathbb{R}^n$.
The following Theorem 5.1 (Theorem 3.2.11. in \cite{rS}) was proved by Shephard \cite{gS} in the case of
a polytope $K$ and by Weil \cite{wW} in the case of compact convex $K$. A strengthening of the theorem 
appeared in Grzybowski, Urba\'nski and Wiernowolski \cite{GUW}.

\vspace{3mm}
\noindent
{\bf Theorem. 5.1 (Shephard--Weil--Schneider criterion).}  
{\it Let $P,K\in \mathcal{B}(\mathbb{R}^n), n\geqslant 2$, 
$P$ be a polytope. Then $P$ is a summand of $K$ if and only 
if the support set $K(u)$ contains a translate of $P(u)$, whenever $P(u)$ 
is an edge of $P$, $u\in S^{n-1}$.} 

\vspace{3mm}
\noindent
The next theorem, a generalization of Theorem 5.1 to an unbounded convex set $K$, 
is based on Schneider's proof from \cite{rS1} presented in Encyclopedia of Mathematics and its Applications 151 \cite{rS}. 
 
\vspace{3mm}
\noindent
{\bf Theorem 5.2 (criterion for a polytopal summand).}  
{\it Let $K\in \mathcal{C}(\mathbb{R}^n), n\geqslant 2$, a recession cone $V$ of $K$ be pointed 
and $P\subset \mathbb{R}^n$ be a polytope. 
Then $P$ is a summand of $K$ if and only 
if every nonempty bounded support set $K(u)$ contains a translate of $P(u)$, whenever $P(u)$ 
is an edge of $P$, $u\in S^{n-1}$.} 

\vspace{3mm}
\noindent
{\it Proof.}  
$\Longrightarrow )$ If a polytope $P$ is a summand of $K$ 
then there exists a set $A\in \mathcal{C}(\mathbb{R}^n)$ 
such that $K=P+A$. If a support set $K(u)$ is nonempty then it is a Minkowski sum 
of respective support sets $K(u)=P(u)+A(u)$. Hence $K(u)$ contains a translate of 
$P(u)$, whether $P(u)$ is an edge or not.

\vspace{3mm}
\noindent
$\Longleftarrow )$
We are going to apply Minkowski duality between convex sets and sublinear functions. 
Basic facts on Minkowski duality are presented in Section 8.
Since the cone $V:=$recc\,$K$ is pointed, the effective domain dom$\,h_K$ has a nonempty interior.
If a difference of support functions $g:=h_K-h_P$ is convex in the interior of dom$\,h_K$ then 
a function $g=h_K-h_P$ is sublinear and lower semicontinuous. Then $K=P+\partial g|_0$, 
and $P$ is a summand of $K$. Hence we need to prove that the function $g$ is convex 
over int\,dom\,$h_K$. Notice that int\,$V^{\circ}\subset$ dom\,$h_K\subset V^{\circ}$, where 
$V^{\circ}$ is a polar of the cone $V$.

\vspace{3mm}
\noindent
Let $x,y\in \textup{int}\,\textup{dom}\,h_K$. 
If $0$ lies between $x$ and $y$ then $0\in \textup{int}\,\textup{dom}\,h_K$ and dom\,$h_K=\mathbb{R}^n$. 
Hence $V=\{0\}$.
This is true only if $K$ is bounded. In this case the polytope $P$ is a summand of $K$ 
by Theorem 5.1.  

\vspace{3mm}
\noindent
Otherwise, lin$\{x,y\}$ is a two-dimensional subspace of $\mathbb{R}^n$. 
Let pr$:\mathbb{R}^n \longrightarrow \textup{lin}\{x,y\}$ be a perpendicular projection. 
Images pr$K$ and pr$P$ of $K$ and $P$ by projection pr are two-dimensional convex sets. 
For any $z\in \textup{lin}\{x,y\}$ equalities 
$h_K(z)=h_{\textup{pr}K}(z)$ and $h_P(z)=h_{\textup{pr}P}(z)$ hold true for respective support functions. 
Assume that every side of the convex polygon pr$P$, that is 
(pr$P)(u)$, $u\in \textup{lin}\{x,y\}$ is equal to an image pr$(P(u))$ of a single 
edge $P(u)$ of the polytope $P$. It simply means that the support set $P(u)$ is an edge of $P$. 
Then if a set (pr$K)(u)=$ pr$(K(u))$, $u\in \textup{lin}\{x,y\}$  is nonempty then (pr$K)(u)$ 
contains a translate of pr$(P(u))=$ (pr$P)(u)$ since $K(u)$ contains a translate of $P(u)$.
Hence by Proposition 4.3 the set pr$P$ is a summand of pr$K$. Thus  
$h_{\textup{pr}K}-h_{\textup{pr}P}$ is a convex function, 
and the function $g=h_K-h_P$ restricted to $\textup{lin}\{x,y\}$ is also convex. Therefore, 
$g(\frac{x+y}{2})\leqslant \frac{g(x)+g(y)}{2}.$ 

\vspace{3mm}
\noindent
If not every side of pr$P$ is equal to a projection of a single 
edge of $P$ then still there exists a sequence $(y_n)$ 
tending to $y$ such that any side of polygon pr$_nP$,  
where pr$_n$ is a perpendicular projection onto the subspace 
$\textup{lin}\{x,y_n\}$,  
is equal to a projection of single edge $P(u)$ of $P$. 

\vspace{3mm}
\noindent
Since the function $g=h_K-h_P$ is continuous in the interior of dom$h_K$, 
we obtain
$$g\left(\frac{x+y}{2}\right)
=\lim\limits_{n\longrightarrow\infty}g\left(\frac{x+y_n}{2}\right)
\leqslant \lim\limits_{n\longrightarrow\infty}\frac{g(x)+g(y_n)}{2}
=\frac{g(x)+g(y)}{2}.$$ 
Since $g$ is continuous in int\,dom\,$h_K$ and $x,y$ are arbitrary, we have 
just proved that $g$ is convex in int\,dom\,$h_K$. 
On the other hand $g=h_K-h_P$ is lower semicontinuous, hence convex in 
all $\mathbb{R}^n$. Therefore, by Theorem 8.1, we obtain $K=P+\partial g|_0$.
\hfill{$\Box$}  

\vspace{3mm}
\noindent 
{\bf Remark 5.3.} Notice that in Theorem 5.2 the assumption of recession cone being 
pointed is necessary. For example let $K$ be a straight line in $\mathbb{R}^n$ and let 
$P$ be any polytope not contained in a straight line parallel to $K$. 
Then $P$ is not a summand of $K$. However, if a support set $K(u)$ is not empty 
then $K(u)=K$, and $K(u)$ is unbounded.  

\vspace{3mm}
\noindent
Let us extend a notion of polytope to unbounded sets sharing a pointed recession cone $V$. 
By $\mathcal{P}_V(\mathbb{R}^n):=\{P+V\,|\,P\in \mathcal{P}(\mathbb{R}^n)\}$, 
where $\mathcal{P}(\mathbb{R}^n)$ is a family of all nonempty polytopes in $\mathbb{R}^n$,
we denote the family of sums of polytopes and the cone $V$.
We call elements of the family $\mathcal{P}_V(\mathbb{R}^n)$ by {\it $V$-polytopes}. 
$V$-polytope is the smallest convex set with a recession cone $V$ containing a given finite set of points.
The following theorem is straightforward corollary from Theorem 5.2.

\vspace{3mm}
\noindent
{\bf Theorem 5.4 (criterion for a $V$-polytopal summand).}  
{\it Let $V$ be a pointed convex cone,  
$K\in \mathcal{C}_V(\mathbb{R}^n), n\geqslant 2$, 
and $P\in \mathcal{P}_V(\mathbb{R}^n)$. 
Then $P$ is a summand of $K$ if and only 
if a nonempty bounded support set $K(u)$ contains a translate of $P(u)$, whenever $P(u)$ 
is an edge of $P$, $u\in S^{n-1}$.} 

\vspace{3mm}
\noindent
Let $A,B\in \mathcal{C}(\mathbb{R}^n)$. We call two bounded support sets $A(u)$ 
and $B(u)$ {\it equiparallel edges} if they are parallel line segments.
Bauer in \cite{cB} gave the following necessary and sufficient criterion for 
reduced pairs of polytopes. 

\vspace{3mm}
\noindent
{\bf Theorem 5.5 (Bauer's criterion for reduced pair of polytopes).}  
{\it A pair $(A,B)$ of polytopes in $\mathbb{R}^n$ is reduced if and only if 
$A$ and $B$ have no equiparallel edges.} 

\vspace{3mm}
\noindent
The next theorem generalizes Bauer's criterion to  
reduced pairs of $V$-polytopes.

%\vspace{3mm}
%\noindent
%{\bf Proposition 4.6.}  
%{\it A pair $(A,B)\in \mathcal{C}_V^2(\mathbb{R}^n)$ is reduced 
%if and only if $A+D=B+C$ implies that $A$ is a summand of $C$ and $B$ is a summand of $D$ 
%for any $C,D\in \mathcal{C}_V(\mathbb{R}^n)$.
%}

%\vspace{3mm}
%\noindent
%%{\bf Proposition 4.7.}  
%{\it Two sets $A,B\in \mathcal{P}_V(\mathbb{R}^n)$ have no 
%equiparallel edges  
%if and only if for all $u\in S^{n-1}$ if $(A+B)(u)$ is an edge 
%then it is a translate of either $A(u)$ or $B(u)$. 
%}

\vspace{3mm}
\noindent
{\bf Theorem 5.6 (criterion for reduced pair of $V$-polytopes).}  
{\it Let $V$ be a pointed convex cone.  
Then a pair $(A,B)\in \mathcal{P}_V^2(\mathbb{R}^n)$ is reduced  
if and only if $A$ and $B$ have no 
equiparallel edges.} 

\vspace{3mm}
\noindent
{\it Proof.}
$\Longleftarrow )$ Let $A$ and $B$ have no equiparallel edges. 
Assume that $A+D=B+C=:E$ for some $C,D\in \mathcal{C}_V(\mathbb{R}^n)$. 
In order to prove that $A+B$ is a summand of $E$, 
let $(A+B)(u)$ be an edge. Since $A$ and $B$ have no equiparallel edges, 
$A(u)$ and $B(u)$ cannot be line segments both  at  the same time. 
Then one of these, say $B(u)$, is a singleton 
and $(A+B)(u)$ is a translate of $A(u)$.  
Hence the set $E(u)=A(u)+D(u)$ contains a translate 
of $(A+B)(u)$. By Theorem 5.4, the set $A+B$ is a summand of $E=A+D=B+C$. 
Therefore, $E=A+B+M$ for some $M\in \mathcal{C}_V(\mathbb{R}^n)$. 
By the cancellation law $(C,D)=(A+M,B+M)$. 

\vspace{3mm}
\noindent 
$\Longrightarrow )$ If $A(u)$ and $B(u)$ are parallel edges then we can construct 
a pair $(A',B')$ equivalent to $(A,B)$ such that $A\subset A', B\subset B'$ and 
no translate of $A(u)$ is contained in $A'(u)$. This construction was given by Bauer in 
Theorem 5.3 \cite{cB} for a pair of polytopes.
\hfill{$\Box$}

\vspace{3mm}
\noindent 
{\bf 6. Application. Minimal representation of a difference of convex functions}

\vspace{3mm}
\noindent
Let $V\subset \mathbb{R}^{n+1}$ be a nontrivial closed convex cone such that 
$V\cap \{x\in \mathbb{R}^{n+1}\,|\, x_{n+1}\geqslant 0\}=\{0\}$.
A pair $(A,B)\in\mathcal{C}^2_V(\mathbb{R}^{n+1})$ is $H$-$minimal$ 
if $(A,B)$ is a minimal element in the family  
$\{(C,D)\in [A,B]\,|\,0\in D\textup{ and }\forall x\in D: x_{n+1}\leqslant 0\}$. 
The definition of $H$-minimality corresponds to Hartman's \cite{pH} definition of a minimal representation  
of a dc-function $f=g-h$, i.e. a difference of convex functions $g$ and $h$, defined on the open unit ball in $\mathbb{R}^n$. 

\vspace{3mm}
\noindent
Let us notice that for two convex and lower semicontinous functions $g,h:\mathbb{R}^n\longrightarrow\mathbb{R}\cup \{+\infty\}$
we can find corresponding cosed convex sets $A,B$ such that
$$g(x)=h_A(x,1), h(x)=h_B(x,1), x\in \mathbb{R}^n.$$ 
The sets $A,B$ are defined by 
\begin{eqnarray*}
A & := & \{(x,t)\in\mathbb{R}^n\times\mathbb{R}\,|\,\forall y:\langle(x,t),(y,1)\rangle\leqslant g(y)\}\\
  &  = & \{(x,t)|\forall y:\langle x,y\rangle+t \leqslant g(y)\},\hspace{45mm} \hfill(\ast \ast\ast)\\
B & := & \{(x,t)|\forall y:\langle x,y\rangle+t \leqslant h(y)\}.\hspace{45mm}
\end{eqnarray*}
Indeed, by Theorem 8.2 the set $A$ is a subdifferential of such a lower semicontinuous sublinear function $\hat{g}$, that 
$\hat{g}(x,t)=tg(x/t), t>0$ and $\hat{g}(x,t)=+\infty, t<0$. 
\noindent
In fact $A=$ hypo$\,(-g^*)$ i.e. the convex set $A$ is equal to a hypograph 
of a function $-g^*$ where $g^*$ is a convex conjugate of $g$ \cite{tR}. We also have 
$B=$ hypo$\,(-h^*)$.

\vspace{3mm}
\noindent
Hartman, defining minimal representation of a dc-function $f=g-h$ 
in section 6 of \cite{pH}, requires that $g,h$ are as small as possible under conditions of
$h\leqslant 0$ and $h(0)=0$.
The function $h$ is non-negative if and only if $0_{\mathbb{R}^{n+1}}\in B$. 
Besides, $h(0)\leqslant 0$ implies $B\subset \mathbb{R}^n\times \mathbb{R}_{-}$. 
If $B\subset \mathbb{R}^n\times \mathbb{R}_{-}$ then 
$h(0)=h_B(0,1)=\sup\limits_{(x,t)\in B}\langle (x,t),(0,1)\rangle
=\sup\limits_{(x,t)\in B}t\leqslant 0$.

\vspace{3mm}
\noindent
Hartman considers dc-function $f=g-h$ defined on an 
interior of a unit ball $\mathbb{B}$ in $\mathbb{R}^n$.
In order to represent convex functions $g,h$ by convex sets 
we extend them outside of int $\mathbb{B}$ by 
$g(x)=h(x)=\infty$ for $x\not\in \mathbb{B}$ and 
$g(x)=\liminf\limits_{y\rightarrow x,\|y\|<1}g(y)$, 
$h(x)=\liminf\limits_{y\rightarrow x,\|y\|<1}h(y)$ for $\|x\|=1$. 

\vspace{3mm}
\noindent
Since effective domains of $g$ and $h$ contain an open Euclidean unit ball 
and are contained in a closed unit ball $\mathbb{B}$, 
the sets $A$ and $B$ share recession cone $V$ defined by 
$V=\{(x,t)\in \mathbb{R}^n\times \mathbb{R}\,|\,t\leqslant -\|x\|_2\}$.
From previous considerations follows the next theorem.

\vspace{3mm}
\noindent
{\bf Theorem 6.1.}
{\it 
A representation $f=g-h$ of a {\rm dc}-function  is minimal according to Hartman
if and only if a pair of sets $(A,B)$, where $A,B$ are defined by $(\ast\ast\ast)$, 
is $H$-minimal in $\mathcal{C}_V^2(\mathbb{R}^{n+1})$ . 
}

\vspace{3mm}
\noindent
If we replace in Hartman's definition an open unit ball with an interior of a closed convex 
set $K$ containing $0$ then corresponding sets $A$ and $B$ share a recession cone $V$ defined by  
$V:=\bigcup\limits_{t\geqslant 0}t(K^{\circ}\times\{-1\})$ 
where $K^{\circ}$ is a polar of $K$. 

%$V=\{(x,t)\in \mathbb{R}^n\times \mathbb{R}\,|\,t\leqslant -\|x\|_2\}
%(=\{(x,t)\in \mathbb{R}^n\times \mathbb{R}\,|\,\frac{-x}{t}\in K^{\circ}\}$.
%$$V=\{(x,t)|\langle(x,t),(\cdot,1)\rangle\leqslant 0\}
%=\{(x,t)|\forall y\in K:\langle(x,t),(y,1)\rangle\leqslant 0\}$$
%$$=\{(x,t)|\forall y\in K:\langle x,y\rangle+t \leqslant 0\},
%=\{(x,t)|\forall y\in K:\langle \frac{x}{-t},y\rangle \leqslant 1\},$$
%$$=\{(x,t)|\sup\limits_{ y\in K}\langle \frac{-x}{t},y\rangle \leqslant 1\}
%=\{(x,t)|h_K(\frac{-x}{t})\leqslant 1\}=\{(x,t)|\frac{-x}{t}\in K^{\circ}\}$$
\vspace{3mm}
\noindent
The following proposition is obvious.

\vspace{3mm}
\noindent
{\bf Proposition 6.2.}
{\it 
A pair $(A,B)\in \mathcal{C}_V^2(\mathbb{R}^{n+1})$ is $H$-minimal if and only if 
it is $0$-minimal and $B\subset \{x\in \mathbb{R}^{n+1}\,|\,x_{n+1}\leqslant 0\}$.
}

\vspace{3mm}
\noindent
In Example 2.10(ii) all equivalent $0$-minimal pairs are $H$-minimal.
Obviously, $0$-minimal pairs may not be $H$-minimal. See the next example.

\vspace{3mm}
\noindent
{\bf Example 6.3.} 
Let $T:\mathbb{R}^3\longrightarrow \mathbb{R}^3$, 
$T(x_1,x_2,x_3):=(x_1,x_2,sx_1+tx_2+x_3), s,t\in \mathbb{R}$. 
Consider convex sets from Example 2.10(ii). 
For $i=0,1,2$ the pairs $(T(\widehat{A}_i),T(\widehat{B}_i))$ 
are $0$-minimal. The pair $(T(\widehat{A}_i),T(\widehat{B}_i))$ is $H$-minimal if and only if 
$-s(p_i)_1-t(p_i)_2\leqslant 0$ and $s(2-(p_i)_1)-t(p_i)_2\leqslant 0$.
For example the pair 
$(T(\widehat{A}_0),T(\widehat{B}_0))$  
where $T(x_1,x_2,x_3):=(x_1,x_2,x_1+x_3)$, $p_0=(0,0)$ 
is $0$-minimal and not $H$-minimal.

\vspace{3mm}
\noindent
Obviously, any pair which is $0$-minimal and not $H$-minimal does not contain 
a $H$-minimal pair.

\vspace{3mm}
\noindent
{\bf Theorem 6.4 (existence of $H$-minimal pairs).}
{\it 
Let $(A,B)\in \mathcal{C}_V^2(\mathbb{R}^{n+1})$. There exists an equivalent $H$-minimal 
pair $(A',B')$ such that $A'\subset A-b, B'\subset B-b$ for some $b\in B$. 
}

\vspace{3mm}
\noindent
{\it Proof.} 
Let $b\in B$ and $b_{n+1}=\max\limits_{x\in B} x_{n+1}$. 
By Theorem 2.1 there exists a $0$-minimal pair $(A',B')$ contained 
in $(A-b,B-b)$. Since $B'\subset B-b\subset \{x\in \mathbb{R}^{n+1}\,|\,x_{n+1}\leqslant 0\}$, 
the pair $(A',B')$ is $H$-minimal.
\hfill{$\Box$}

\vspace{3mm}
\noindent
{\bf Remark 6.5.}
It is possible that among equivalent pairs of sets a $H$-minimal pair is unique even if 
this pair does not have the property of translation. For example 
the pair $(T(\widehat{B}_0)-(1,1,1),T(\widehat{A}_0)-(1,1,1))$ from Example 6.3, 
where $T(x_1,x_2,x_3):=(x_1,x_2,x_2+x_3)$, $p_0=(0,0)$,
is a unique $H$-minimal pair in the quotient class $[T(\widehat{B}_0),T(\widehat{A}_0)]$. 
Notice that $T(\widehat{A}_0)-(1,1,1)=$\,conv$\{(0,0,0),(-1,-1,-1),(1,-1,-1)\}+V$, 
$T(\widehat{B}_0)-(1,1,1)=$\,conv$\{(-1,-1,-1),(1,-1,-1)\}+V$. Convex functions corresponding 
to these two sets are $g(x_1,x_2):=|x_1|-x_2-1$ and $h(x_1,x_2):=\max(0,|x_1|-x_2-1)$. 
They are the unique Hartman-minimal 
convex functions, such that $f(x_1,x_2):=\min(0,|x_1|-x_2-1)=g(x_1,x_2)-h(x_1,x_2)$. 

\vspace{3mm}
\noindent
{\bf Proposition 6.6.}
{\it 
Let a pair $(A,B)\in \mathcal{C}_V^2(\mathbb{R}^{n+1})$ be reduced and  
$V\cap \{x\in \mathbb{R}^{n+1}\,|\,x_{n+1}\geqslant 0\}=L_V$.
Then a pair $(A-x,B-x), x\in B$ is $H$-minimal if and only if   
$x_{n+1}=\sup\limits_{y\in B}y_{n+1}=h_B(u)$, where $u=(0,...,0,1)\in \mathbb{R}^{n+1}$.
}

\vspace{3mm}
\noindent
{\it Proof.} 
By Theorem 3.3 the pair $(A,B)$ has property of translation. Proposition follows 
from criterion of $0$-minimality in Proposition 2.6(b) and from
characterization of $H$-minimality in Proposition 6.2.
\hfill{$\Box$}

\vspace{3mm}
\noindent
The following example shows, that reducibility (property of translation) in the assumptions of 
Proposition 2.6(b) is essential. 

\vspace{3mm}
\noindent
{\bf Example 6.7.} Consider a pair $(\widehat{A}_1,\widehat{B}_1)$ in Example 2.10(ii). 
This pair is $H$-minimal but not reduced. 
If $x\in \widehat{B}_1$ then $(\widehat{A}_1-x,\widehat{B}_1-x)$ is $H$-minimal if and only if $x=0$. 
Still $(\widehat{A}_0-x,\widehat{B}_0-x)$ is $H$-minimal if and only if 
$x\in (\widehat{B}_0)_{*}$, i.e. 
$x=(x_1,...,x_{n+1})\in \widehat{B}_0$ and $x_{n+1}=\sup\limits_{y\in \widehat{B}_0}y_{n+1}=0$.

\vspace{3mm}
\noindent
{\bf 7. Appendix. Minimal pairs of closed bounded convex sets}

\vspace{3mm}
\noindent
Let $\mathcal{B}(\mathbb{R}^n)$ be a family of all nonempty compact convex sets, i.e. convex bodies.
The idea of treating compact convex sets as numbers or, rather, as vectors goes back to 
Minkowski \cite{hM}.  
A semigroup of nonempty bounded closed convex subsets $\mathcal{B}(X)$ 
of a vector space $X$ was embedded into a topological vector space in the case of a normed 
space $X$ by R\aa dstr\"om \cite{hR}, a locally convex space $X$ by H\"ormander 
\cite{lH} and a topological vector space by Urba\'nski \cite{rU}.
The embedding was possible thanks to an order cancellation law: 
$$A+B\subset \textup{\,cl\,}(B+C) \Longrightarrow A\subset C 
\textup{\,\,\, for \,\,\,} A,B,C\in \mathcal{B}(X).$$ 
For a concise proof of an order cancellation law in a more general setting 
we refer the reader to Proposition 5.1 in \cite{GKKU}. 
Convex sets are embedded into Minkowski--R\aa dstr\"om--H\"ormander space 
$\widetilde{X}=\mathcal{B}^2(X)/_{\sim}$ of quotient classes, where 
a relation of equivalence is defined by 
$(A,B)\sim(C,D):\Longleftrightarrow \textup{\,cl\,}(A+D)=\textup{\,cl\,}(B+C)$. 

\vspace{3mm}
\noindent
A new motivation to study pairs of convex sets came from quasidifferential calculus 
of Demyanov and Rubinov \cite{DR1,DR2}, 
where a quasidifferential $Df(x_0)$ is a pair of convex sets 
$(A,B)=(\underline{\partial}f|_{x_0},\overline{\partial}f|_{x_0})$ 
called sub- and superdifferential. Rather than a pair of sets $(A,B)$ 
a quasidifferential 
is a quotient class $[A,B]:=[(A,B)]_{\sim}$. 

\vspace{3mm}
\noindent
The best representation of a quotient class $[A,B]$ is a {\it reduced 
pair}, i.e. a pair $(A,B)$ such that $[A,B]=\{(A+C,B+C)\,|\,C\in \mathcal{B}(X)\}$. 
Then all translates of $(A,B)$ give all minimal elements of $[A,B]$.  
Reduced pairs were studied by Bauer \cite{cB}. However, not every quotient class $[A,B]$ 
contains a reduced pair. We say that a pair $(A,B)$, or a quotient class $[A,B]$ 
has {\it property of translation} if all minimal pairs in $[A,B]$ 
are translates of each other. There exist not reduced minimal pairs that 
have property of translation.   
The following theorem holds true. 

\vspace{3mm}
\noindent
{\bf Theorem 7.1.} (\cite{GUa,PSU})
{\it 
Let $X$ be a reflexive Banach space. For every pair $(A,B)\in \mathcal{B}^2(X)$, 
there exists an inclusion-minimal pair $(C,D)\in [A,B]$ such that $C\subset A, D\subset B$.} 

\vspace{3mm}
\noindent
Caprari and Penot \cite{CP} proved existence of inclusion minimal pairs in a quotient class 
$[A,B]\in \mathcal{C}(X)\times \mathcal{K}(X)/_{\sim}$ where $\mathcal{K}(X)$ is a family 
of all nonempty compact convex subsets of a locally convex vector space $X$.

\vspace{3mm}
\noindent
{\bf Theorem 7.2.} (\cite{cB,jG,sS})
{\it 
Let $(A,B)\in \mathcal{B}^2(\mathbb{R}^2)$. A minimal pair 
in $[A,B]$ is unique up to translation.} 

\vspace{3mm}
\noindent
Theorem 7.2 basically states that every minimal pair of two-dimensional compact convex sets 
has property of translation.  

\vspace{3mm}
\noindent
{\bf Example 7.3.} (\cite{jG})
In $\mathbb{R}^3$ we have equivalent minimal pairs not connected by translation .

\begin{displaymath}
\begin{array}{ccccc}

\setlength{\unitlength}{.4mm}
\begin{picture}(50,60)(0,0)
\thinlines{
\put(16,5){\makebox(0,0){{\bf $A$}}}
\put(00,18){\line(0,1){12}}
\put(20,18){\line(0,1){12}}
\put(10,0){\line(0,1){12}}
\qbezier(00,30)(00,30)(10,48)
\qbezier(00,30)(00,30)(10,12)
\qbezier(00,18)(00,18)(10,0)
\qbezier(20,30)(20,30)(10,48)
\qbezier(20,30)(20,30)(10,12)
\qbezier(20,18)(20,18)(10,0)

\put(46,09){\makebox(0,0){{\bf $B$}}}
\put(30,18){\line(0,1){12}}
\put(50,18){\line(0,1){12}}
\qbezier(30,30)(30,30)(40,12)
\qbezier(50,30)(50,30)(40,12)
\qbezier(30,30)(30,30)(40,36)
\qbezier(30,18)(30,18)(40,12)
\qbezier(50,30)(50,30)(40,36)
\qbezier(50,18)(50,18)(40,12)
\put(30,30){\line(1,0){20}}
}
\end{picture}
 
&\,&

\setlength{\unitlength}{.4mm}
\begin{picture}(50,60)(0,0)
\thinlines{
\put(18,5){\makebox(0,0){{\bf $C$}}}
\put(00,18){\line(0,1){28}}
\put(20,18){\line(0,1){28}}
\put(10,0){\line(0,1){12}}
\qbezier(00,46)(00,46)(10,64)
\qbezier(00,46)(00,46)(10,12)
\qbezier(00,18)(00,18)(10,0)
\qbezier(20,46)(20,46)(10,64)
\qbezier(20,46)(20,46)(10,12)
\qbezier(20,18)(20,18)(10,0)
\put(00,18){\line(5,-3){10}}
\put(20,18){\line(-5,-3){10}}
\put(00,46){\line(1,0){20}}

\put(47,10){\makebox(0,0){{\bf $D$}}}
\put(30,18){\line(0,1){28}}
\put(50,18){\line(0,1){28}}
\put(30,46){\line(1,0){20}}
\qbezier(30,46)(30,46)(40,12)
\qbezier(50,46)(50,46)(40,12)
\qbezier(30,46)(30,46)(40,52)
\qbezier(30,18)(30,18)(40,12)
\qbezier(50,46)(50,46)(40,52)
\qbezier(50,18)(50,18)(40,12)
}
\end{picture}

&\,&

\setlength{\unitlength}{.4mm}
\begin{picture}(90,36)(0,0)
\thinlines{
\put(38,5){\makebox(0,0){{\bf $E$}}}
\put(10,0){\line(0,1){12}}
\put(30,0){\line(0,1){12}}
\put(10,12){\line(1,0){20}}
\put(10,0){\line(1,0){20}}
\put(10,36){\line(1,0){20}}
\qbezier(30,0)(30,0)(40,18)
\qbezier(10,36)(10,36)(00,18)
\qbezier(20,30)(20,30)(10,12)
\qbezier(20,30)(20,30)(30,36)
\qbezier(00,18)(00,18)(10,12)
\qbezier(20,30)(20,30)(10,36)
\qbezier(40,18)(40,18)(30,12)
\qbezier(00,18)(00,18)(10,0)
\qbezier(20,30)(20,30)(30,12)
\qbezier(30,36)(30,36)(40,18)

\put(86,10){\makebox(0,0){{\bf $F$}}}
\put(60,12){\line(1,0){20}}
\put(60,24){\line(1,0){20}}
\qbezier(50,18)(50,18)(60,24)
\qbezier(50,18)(50,18)(60,12)
\qbezier(90,18)(90,18)(80,24)
\qbezier(90,18)(90,18)(80,12)

}
\end{picture}
\end{array}
\end{displaymath}

\begin{center}
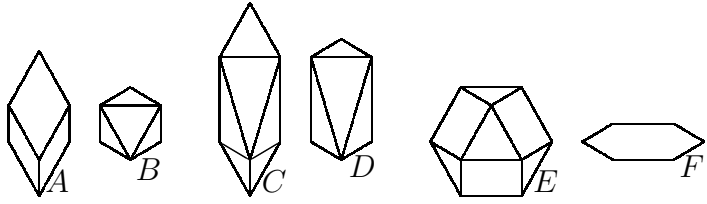
Figure 7.1. Three equivalent minimal pairs not connected by translation.
\end{center}

\vspace{3mm}
\noindent
In Figure 7.1 the solid $B$ is a regular octahedron, $D$ is an elongated octahedron, 
$F$ is a hexagon, $A$ is a rhombohedron and $E$ is a cuboctahedron.
 
\vspace{3mm}
\noindent
In \cite{GKU,GPU,PU} more quotient classes $[A,B]$ with no unique minimal pair 
were found in $\mathbb{R}^3$. However, the set of all equivalent minimal pairs 
was never effectively described for a quotient class $[A,B]$ with no unique minimal 
element. All these results enabled calculus of pairs of convex sets in a way 
analogous to fractional arithmetics \cite{PUB}. 

\vspace{3mm}
\noindent
The following theorem states a necessary and sufficient criterion for minimal pairs of 
convex polygons.

\vspace{3mm}
\noindent
{\bf Theorem 7.4 (Theorem 3.5 in \cite{GPU2009}).}
{\it A pair $(A,B)$ of flat polytopes is minimal
if and only if $A$ and $B$ have at most one pair of parallel edges that lie on the same side
of the polytopes.}

\vspace{3mm}
\noindent
{\bf 8. Appendix. Minkowski duality}

\vspace{3mm}
\noindent
In this section we present Minkowski duality between closed convex sets and 
sublinear functions.

\vspace{3mm}
\noindent
A 4-tuple $(X,\mathbb{R}_{+},+,\cdot)$, where an operation of addition '$+$'  
and of multiplication by nonnegative numbers '$\cdot$' are defined 
for elements of the set $X$, is called an {\it abstract convex cone} 
if (1) the pair $(X,+)$ is a commutative group and for all 
$x,y\in X$ and all $s,t\geqslant 0$ we have (2) $1x=x$, (3) $0x=0$, 
(4) $s(tx)=(st)x$, (5) $t(x+y)=tx+ty$ and (6) $(s+t)x=sx+tx$.  
 
\vspace{3mm}
\noindent
If a set $A$ belongs to the family $\mathcal{C}(\mathbb{R}^n)$ of all 
nonempty closed convex subsets of $\mathbb{R}^n$ then 
its {\it support function} $h_A$ is defined by 
$$h_A:=\sup\limits_{a\in A}\langle a,\cdot\rangle.$$
Minkowski addition $A\dot{+}B$ of sets belonging to $\mathcal{C}(\mathbb{R}^n)$ 
is defined by $A\dot{+}B:=\textup{cl}(A+B)$. 
Obviously, if one of these sets is bounded then $A\dot{+}B=A+B$.

\vspace{3mm}
\noindent
If a function $h:\mathbb{R}^n\longrightarrow \mathbb{R}\cup \{\infty\}$ belongs 
to the family $\mathcal{Sub}^{\infty}_{lsc}(\mathbb{R}^n)$ of all 
sublinear (positively homogenous and convex) lower semicontinuous functions 
then its {\it subdifferential} at $0$ is a closed convex set 
defined by $\partial h|_0:=\{x\in\mathbb{R}^n\,|\,\langle x,\cdot\rangle\leqslant h\}$.

\vspace{3mm}
\noindent
Both 4-tuples $(\mathcal{C}(\mathbb{R}^n),\mathbb{R}_{+},\dot{+},\cdot)$ and 
$(\mathcal{Sub}^{\infty}_{lsc}(\mathbb{R}^n),\mathbb{R}_{+},+,\cdot)$ are 
abstract convex cones and Minkowski duality establishes isomorphic 
relationship between these cones.

\vspace{3mm}
\noindent
{\bf Theorem 8.1.} {\it 
The mapping $\mathcal{Sub}^{\infty}_{lsc}(\mathbb{R}^n)\ni h 
\longmapsto \partial h|_0\in \mathcal{C}(\mathbb{R}^n)$ 
is an isomorphic bijection from an abstract convex cone
$\mathcal{Sub}^{\infty}_{lsc}(\mathbb{R}^n)$ onto an abstract convex cone 
$\mathcal{C}(\mathbb{R}^n)$. 
The mapping $\mathcal{C}(\mathbb{R}^n)\ni A 
\longmapsto h_A\in\mathcal{Sub}^{\infty}_{lsc}(\mathbb{R}^n)$ is an inverse mapping.
Moreover, a restriction of the mapping to the subfamily of finite sublinear functions 
$\mathcal{Sub}(\mathbb{R}^n)$ is an isomorphic bijection from a subcone $\mathcal{Sub}(\mathbb{R}^n)$
onto a subcone $\mathcal{B}(\mathbb{R}^n)$ of all nonempty bounded closed convex sets.
}

\vspace{3mm}
\noindent
Theorem 8.1 is stated in \cite{GPU70} in a general case for a dual pair $(X, Y )$ of 
linear spaces over $\mathbb{R}$ where 
$\langle\cdot,\cdot\rangle$ is such a bilinear function, 
that functions $\{\langle y, \cdot\rangle \}_{y\in Y}$ separate points 
in $X$ and functions $\{\langle\cdot,x\rangle\}_{x\in X}$ separate points in $Y$.

\vspace{3mm}
\noindent
Let $h:\mathbb{R}^n\longrightarrow\mathbb{R}\cup\{\infty\}$ be a sublinear lower semicontinous 
function. An effective domain dom\,$h\subset \mathbb{R}^n$ is a convex cone, its 
closure cl\,dom\,$h$ is a closed convex cone.  

\vspace{3mm}
\noindent
Let $V$ be a closed convex cone in $\mathbb{R}^n$. 
A characteristic function $\chi_V$ is defined by 
$%\begin{displaymath}
\chi_V(x):=\left\{
\begin{array}{cc}
0,                                       & x\in V,\\
\infty,                                  & x\not\in V. 
\end{array}
\right.
$%\end{displaymath}
A subdifferential of $\chi_V$ at $0$ coincides with a polar cone $V^{\circ}$.
Therefore, $\partial(\chi_V)|_0=V^{\circ}$, and $h_V=\chi_{V^\circ}$.

\vspace{3mm}
\noindent
By $\mathcal{Sub}^{\infty}_{lsc,V}(\mathbb{R}^n)$ we denote a subfamily 
of sublinear functions with finite values in the relative interior of $V$ and
infinite outside of $V$. Values of such a function on the relative boundary of 
$V$ are determined by its values in the relative interior.  
A 4-tuple $(\mathcal{Sub}^{\infty}_{lsc,V}(\mathbb{R}^n),\mathbb{R}_{+},+,\cdot)$ 
is an abstract convex cone after modifying multiplication by 0 in the following way 
$0h=:\chi_V$.

\vspace{3mm}
\noindent
By $\mathcal{C}_{V}(\mathbb{R}^n)$ we denote all closed convex sets $A$ having 
their recession cone recc$A:=A\dot{-}A=\{x\,|\,x+A\subset A\}$ equal to $V$. 
Again the family $\mathcal{C}_{V}(\mathbb{R}^n)$ is an abstract convex cone 
after modifying multiplication by $0$ with a formula $0A=V$ \cite{sR}.

\vspace{3mm}
\noindent
{\bf Theorem 8.2.} {\it 
The mapping $\mathcal{Sub}^{\infty}_{lsc,V}(\mathbb{R}^n)\ni h 
\longmapsto \partial h|_0\in \mathcal{C}_{V^{\circ}}(\mathbb{R}^n)$ 
is an isomorphic bijection from an abstract convex cone
$\mathcal{Sub}^{\infty}_{lsc,V}(\mathbb{R}^n)$ onto an abstract convex cone 
$\mathcal{C}_{V^{\circ}}(\mathbb{R}^n)$. 
The mapping $\mathcal{C}_{V^{\circ}}(\mathbb{R}^n)\ni A 
\longmapsto h_A\in\mathcal{Sub}^{\infty}_{lsc,V}(\mathbb{R}^n)$ an is an inverse mapping.
}

\vspace{3mm}
\noindent
{\it Proof.} In view of Theorem 8.1 it is enough to prove that for any function 
$h\in \mathcal{Sub}^{\infty}_{lsc}(\mathbb{R}^n)$ closed 
convex cones $V_1=$cl\,dom$\,h$ and $V_2=$recc$(\partial h|_0)$ are mutually polar.  

\vspace{3mm}
\noindent
Notice that $h=h+\chi_{V_1}$. By Theorem 8.1 we obtain $\partial h|_0=
\partial h|_0+\partial (\chi_{V_1})|_0=\partial h|_0+V_1^{\circ}$.
Then $V_1^{\circ}\subset\partial h|_0\dot{-}\partial h|_0=$recc$(\partial h|_0)=V_2$.

\vspace{3mm}
\noindent
On the other hand $\partial h|_0=\partial h|_0+V_2$. By Theorem 8.1 we get 
$h=h+h_{V_2}=h+\chi_{V_2^{\circ}}$. Hence dom$\,h\subset V_2^{\circ}$, 
and $V_1\subset V_2^{\circ}$. 
Therefore, $V_2\subset V_1^{\circ}$. 
\hfill{$\Box$}

\end{document}